\documentclass[11pt,leqno]{amsart}
\usepackage{amscd}
\usepackage{amssymb}
\usepackage{amsfonts}
\usepackage{latexsym}
\usepackage{verbatim}
\usepackage{amsthm}
\numberwithin{equation}{section}

\theoremstyle{plain}
\newtheorem{theorem}{Theorem}[section]

\begin{document}
\title{Finiteness of Ulam polynomials}
\author{Antonio J. Di Scala and \'Oscar Maci\'a}
\date{July, 2008}
\address{Dipartimento di Matematica
\\ Politecnico di Torino \\ Corso Duca degli Abruzzi 24, 10129 Torino,
Italy} \email{antonio.discala@polito.it}
\email{oscarmacia@calvino.polito.it}
\subjclass[2000]{Primary 11C08 ; Secondary 12D10}
\thanks{The work of O.M. was supported by the Spanish Ministry of Science and Education through a postdoctoral fellowship.}
\begin{abstract}
A polynomial whose coefficients are equal to its roots is called a Ulam polynomial.
In this paper we show that for a given degree $n$ there exists a finite number of Ulam polynomials of degree $n$.
\end{abstract}
\maketitle
\newcommand\C{{\mathbb C}}
\newcommand\R{{\mathbb R}}
\newcommand\Z{{\mathbb Z}}
\newcommand\T{{\mathbb T}}
\newcommand\GL{{\rm GL}}
\newcommand\SL{{\rm SL}}
\newcommand\SO{{\rm SO}}
\newcommand\Sp{{\rm Sp}}
\newcommand\U{{\rm U}}
\newcommand\SU{{\rm SU}}
\newcommand{\Gdue}{{\rm G}_2}
\newcommand\re{\,{\rm Re}\,}
\newcommand\im{\,{\rm Im}\,}
\newcommand\id{\,{\rm id}\,}
\newcommand\tr{\,{\rm tr}\,}
\renewcommand\span{\,{\rm span}\,}
\newcommand\Ann{\,{\rm Ann}\,}
\newcommand\Hol{{\rm Hol}}
\newcommand\Ric{{\rm Ric}}
\newcommand\nc{\widetilde{\nabla}}
\renewcommand\d{{\partial}}
\newcommand\dbar{{\bar{\partial}}}
\newcommand\s{{\sigma}}
\newcommand\sd{{\bigstar_2}}
\newcommand\K{\mathbb{K}}
\renewcommand\P{\mathbb{P}}
\newcommand\D{\mathbb{D}}
\newcommand\al{\alpha}
\newcommand\f{{\varphi}}
\newcommand\g{{\frak{g}}}
\renewcommand\k{{\kappa}}
\renewcommand\l{{\lambda}}
\newcommand\m{{\mu}}
\renewcommand\O{{\Omega}}
\renewcommand\t{{\theta}}
\newcommand\ebar{{\bar{\varepsilon}}}
\newcommand{\dis}{\di^{\circledast}}
\newcommand{\fo}{\mathcal{F}}
\newcommand{\op}{\overline{\parzial}}
\newcommand{\wn}{\widetilde{\nabla}}
\section{Introduction}

In \cite[p. 31]{Ulam} S. Ulam proposed to study some properties of the following map $T_n : \C^n \rightarrow \C^n$
\[ T_n : x'_j = (-1)^j \sigma_j(x_1,\cdots,x_n) \]
where $\sigma_j$ is the $j$th elementary symmetric function. Hence, $x'_1, \cdots,x'_n$ are the coefficients of the equation \[ z^n + x'_1 z^{n-1}+ \cdots + x'_{n-1}z + x'_n = (z-x_1)(z - x_2) \cdots (z - x_{n-1}) (z-x_n) = 0 \, \]
whose roots are $x_1,\cdots,x_n$. He wrote:
\begin{quote}
M\emph{any of the  statements about algebraic equations are translatable into the elementary properties of this mapping. Thus, Gauss' theorem on the existence of roots is simply the statement that $T_n$ is a mapping (many-one) on $E^n$ to all of $E^n,$ where $E$ is the complex plane. The points ``constructible by ruler and compass'' are related to those resulting from iteration of the inverse transformation $T^{-1}_n$ where $n=2.$\\ However, the topological nature of this transformation does not seem to have been very thoroughly investigated. For example, what are the nontrivial fixed points $p=T_n(p)$? The origin is always a fixed point, but there are others, e.g., $p=(1,-2)$ when $n=2.$} \end{quote}

A monic polynomial $P \in \C[z]$ of degree $n$ is called a Ulam polynomial if its roots are a fixed point of the map $T_n$. Namely, $P(z) = (z - x_1) \cdots (z - x_n)$ is Ulam if $T_n(x_1,\cdots,x_n) = (x_1,\cdots,x_n)$. Let us call $U_n$ the set of Ulam polynomials of degree $n$. We will regard $U_n$ as a subset of  $\C^n$.

In  \cite{Stein}  Stein  showed that $U_n$ contains no nontrivial ($P(0) \neq 0$) real polynomial for $n\geq 5$.  Here is our main result.

\begin{theorem}\label{main}
$U_n$ is finite for all $n$.
\end{theorem}

We recall the obvious fact that the slightly different problem arising from considering polynomials whose roots are equal to the opposite of their coefficients allows the possibility of having infinite families of polynomials satisfying this condition for any degree e.g., $P(z) = z - \alpha$.

\section{Proof of Theorem \ref{main}}

Let $V(f) \subset \C^{n+1}$ be the hypersurface defined by $f(z,x_1,x_2, \cdots, x_n) = z^n + x_1 z^{n-1}+ \cdots + x_{n-1}z + x_n-(z-x_1)(z - x_2) \cdots (z - x_{n-1}) (z-x_n)$. Notice that if $P \in U_n$ then $P$ produces a line $\mathcal{L}_P$ in $V(f)$. Indeed, if $P(z)= (z - x_1)\cdots(z-x_n)$ then the line $z \mapsto (z,x_1,x_2,\cdots,x_n)$ is contained in $V(f)$.
We can regard $\C^{n+1}$ as an affine chart of the projective space $\mathbb{P}^{n+1}$ i.e., in homogeneous coordinates $[z:x_1:x_2:\cdots:x_n : 1]$. Assume now that $U_n$ is infinite. Since $U_n$ is a Zariski closed set, there exists an unbounded sequence of point $P_j \in U_n$.\\
Then there exists a sequence $\mathcal{L}_j$ of lines in $V(f)$ associated to a sequence of Ulam polynomials $P_j$ for $j=1,2,\cdots$. Since $G(1,n+1)$, the Grassmannian of lines of $\mathbb{P}^{n+1}$, is compact we can assume that $\mathcal{L}_j$ is convergent to a limit $\mathcal{L}$.
Notice that $\mathcal{L}\subset \P^{n+1}=\C^{n+1}\cup\P^n.$ We claim that $\mathcal{L}\cap\C^{n+1}=\emptyset.$ Indeed, if $\mathcal{L}\cap\C^{n+1}\neq\emptyset,$ then $\mathcal{L}\subset \C^{n+1},$ and this contradicts the fact that $\{P_j\}$ is unbounded. Thus, $\mathcal{L}\subset \P^n.$\\
Let us compute $\mathcal{L}$ in homogeneous coordinates. The homogeneous equation
\[ z^n + x_1 z^{n-1}+ \cdots + x_{n-1}zw^{n-2} + x_nw^{n-1} =(z-x_1)(z - x_2) \cdots (z - x_{n-1}) (z-x_n)\]
implies that a point in $\mathcal{L}$ with $w=0$ satisfies
\[ z^n + x_1 z^{n-1} = z^{n-1}(z + x_1) = (z-x_1)(z - x_2) \cdots (z - x_{n-1}) (z-x_n) \, \, .\]
Then $x_1 = x_2 = \cdots = x_n = 0$. Notice that $z$ is still a parameter for the limit line. Thus an arbitrary point $p \in \mathcal{L}$ has homogeneous coordinates $[z:0:0:0:\cdots:0] = [1 : 0 : 0 : \cdots : 0]$. This shows that the limit $\mathcal{L}$ consists of a single point which gives a contradiction since $\mathcal{L}$ is a line. $\Box$

\section{Acknowledgements}
The authors would like to thank Martin Sombra and Carlos D'Andrea.

\end{document}